\documentclass{amsart}[12pt]
\usepackage{a4}
\usepackage[T1]{fontenc}
\usepackage[latin1]{inputenc}
\usepackage{hyperref}
\usepackage[dvips]{graphics}
\usepackage{amsmath}
\usepackage{amssymb}
\usepackage{amsopn}
\usepackage{amsbsy}
\usepackage{amstext}
\usepackage{latexsym}
\usepackage{amsfonts}
\usepackage{array}
\usepackage{epsfig}
\usepackage{mathrsfs}
\title[About the fractional parts of the powers of rational numbers]{About the fractional parts of the powers of the rational numbers}
\author[B. FARHI]{Bakir FARHI}
\date{\bf November 20\textsuperscript{th}, 2006}

\newtheorem{thm}{Theorem}

\let\epsilon=\varepsilon
\def\EMts{\mspace{.3mu}}

\def\norm#1{{\left\Vert{\EMts\EMts #1 \EMts\EMts}\right\Vert}}
\def\EMdash{\leavevmode\hbox to 7.5mm{\vrule height .63ex depth -.59ex
    width 5.4mm\hfill}}

\begin{document}
\maketitle \baselineskip=6mm
\begin{center}
Département de Mathématiques, Université du Maine, \\
Avenue Olivier Messiaen, 72085 Le Mans Cedex 9, France. \\
Bakir.Farhi@univ-lemans.fr
\end{center}
\begin{abstract}
Let $p / q$ $(p , q \in \mathbb{N}^*)$ be a positive rational number
such that $p > q^2$. We show that for any $\epsilon
> 0$, there exists a set $A(\epsilon) \subset [0 , 1[$, with finite border and with Lebesgue measure $< \epsilon$, for which the set of
positive real numbers $\lambda$ satisfying $<\lambda (p / q)^n> \in
A(\epsilon)$ $(\forall n \in \mathbb{N})$ is uncountable.
\end{abstract}
\hfill{--------} 
\section{Introduction}
Throughout this paper, $\lfloor x \rfloor$ and $<x>$ denote
respectively the integer and the fractional part of a given real
number $x$ and $\norm{x}$ denotes the distance from $x$ to the
nearest integer.

The distribution of the fractional parts of rational numbers is an
important and fertile subject by the numerous problems connected
to which. The most famous of those problems is the determination
of the exact values $g(n)$ $(n \in \mathbb{N}^*)$ of
Hilbert-Waring's theorem\footnote{Here $g(n)$ denotes the least
positive integer such that any positive integer can be written as
a sum of $g(n)$ $n$\textsuperscript{th} powers of non-negative
integers.}. This last problem is oddly related to the distribution
of the sequence ${\{<(3 / 2)^n>\}}_{n \geq 1}$. Indeed, we know
for instance (see e.g. \cite{hua2}) that if an integer $n
> 6$ satisfies:
\begin{equation}
<(3 / 2)^n> \leq 1 - (1 / 2)^n \left(\lfloor(3 / 2)^n\rfloor +
3\right) , \tag{$*$}
\end{equation}
then we have:
\begin{equation}
g(n) = 2^n + \lfloor(3 / 2)^n\rfloor - 2 . \tag{$**$}
\end{equation}
A variant of this last result consists to replace $(*)$ by:
\begin{equation*}
\norm{(3 / 2)^n} > (3 / 4)^n
\end{equation*}
and always conclude to $(**)$ (see also \cite{ben} for a survey of results on Waring's problem).\\
Further, K. Mahler \cite{mah1} showed that $(*)$ is satisfied for
any sufficiently large $n$, i.e., for any integer greater than
some positive integer $n_0$\footnote{It must be noted that none
estimate from above of $n_0$ is known. This stems from the
infectivity of Roth's diophantine theorem used by K. Mahler.}. It
follows that $(**)$ holds for any sufficiently large integer $n$.
The problem of knowing starting from which positive integer
precisely $(**)$ is satisfied is always open. Since it is well
known that $(**)$ holds for the first small values of $n$, it is
very probable that it holds for all $n \in \mathbb{N}^*$.

One of the first important studies about the distribution of the
fractional parts of the powers of rational numbers was that of T.
Vijayaraghavan \cite{vij}. This last showed that
\begin{quote}
``for any rational number $p / q$ $(p , q \in \mathbb{N}^* , p
> q)$, the sequence ${\{<(p / q)^n>\}}_{n \in \mathbb{N}}$ has an infinitely many limit points.''
\end{quote}
It is however curious to note that Vijayaraghavan's method does
not make it possible to say if there is an infinitely many limit
points of the sequence ${\{<(3 / 2)^n>\}}_n$ in each interval $[0
, 1 / 2[$ and $[1 / 2 , 1[$. In \cite{vij}, one conjectures that
this is the case and one also conjectures that one has:
\begin{equation}
\limsup_{n \rightarrow + \infty} <(3 / 2)^n> - \liminf_{n
\rightarrow + \infty} <(3 / 2)^n> > \frac{1}{2} . \tag{$***$}
\end{equation}
K. Mahler \cite{mah2} formulated the conjecture which stating that
\begin{quote}
``There is not real positive number $\lambda$ such that:
$$<\lambda (3 / 2)^n> < \frac{1}{2} ~~~~~~ (\forall n \in \mathbb{N}) .\text{''}$$
\end{quote}
Such numbers $\lambda$ are called ``Z-nombres''. The
conjecture of K. Mahler thus amounts to say that Z-numbers do not exist.\\
In the direction to prove $(***)$, L. Flatto, J. C. Lagarias and A.
D. Pollington \cite{flp} showed that for any real positive number
$\lambda$ and any rational number $p / q$ (with $p , q \in
\mathbb{N}^* , p > q$), we have:
$$\limsup_{n \rightarrow + \infty} <\lambda (p / q)^n> - \liminf_{n \rightarrow + \infty} <\lambda (p / q)^n>
\geq \frac{1}{p} .$$ This result is generalized very recently by
A. Dubikas \cite{dub} which deal with algebraic numbers instead of
rational numbers. One of the important consequences of the results
of \cite{dub} is that the sequence ${\{<\lambda (3 / 2)^n>\}}_n$
$(\lambda
> 0)$ admits at least one limit point in the interval
$[0.238117\dots , 0.761882\dots]$. \\
Further, A. D. Pollington \cite{pol} had already shown that the
set of real positive numbers $\lambda$ satisfying $\norm{\lambda
(3 / 2)^n} > \frac{4}{65}$ $(\forall n \in \mathbb{N})$ has
positive Hausdorff dimension, consequently it is uncountable.

In \cite{far}, the author showed some results concerning the
distribution modulo $1$ of geometric progressions. The basic
result states that for all $r > 1$, the set
$$E_r := \left\{\lambda > 0 ~|~ \norm{\lambda r^n} \leq \frac{1}{r - 1} ~~ (\forall n \in \mathbb{N})\right\}$$
is uncountable.\\
One of the consequences of this result is that for all $r \in ]3 ,
+ \infty[ \cup \{2 , 3 , 5 / 2\}$, the set of real positive
numbers $\lambda$ for which the sequence $\{\lambda r^n\}_{n \in
\mathbb{N}}$ is not dense modulo $1$, is uncountable (cf.
Corollary 2 or 2.2 of \cite{far}). We remark here that the
rational numbers $2 , 3$ et $5 / 2$ appear as particular cases! If
we look closely at the method used to deal with the case $r = 5 /
2$, we see that the secrecy of this particularity is that these
three numbers all have the form $p / q$ $(p , q \in \mathbb{N}^*)$
with ${\bf p > q^2}$. In this note, we take up the idea of the
proof of Corollary 2.2 (case $r = 5/2$) of \cite{far} and apply it
to the more general case of the rational numbers $p / q$ $(p , q
\in \mathbb{N}^* , p > q^2)$. We will establish in this way a new
information concerning the distribution modulo $1$ of geometric
progressions which have for common ratio such rational numbers.
\section{The Result}
\begin{thm}\label{t1}
Let $r = p / q$ $(p , q \in \mathbb{N}^*)$ be a rational positive
number such that $p > q^2$. Then, for any $\epsilon
> 0$, there exists a set $A(\epsilon) \subset [0 , 1[$, with finite border and with Lebesgue measure $< \epsilon$, for which the set of
positive real numbers $\lambda$ satisfying:
$$<\lambda r^n> \in A(\epsilon) \text{\hspace{1cm}} (\forall n \in \mathbb{N})$$
is uncountable.
\end{thm}\noindent
{\bf Proof.---} The fact $p > q^2$ implies
$$\lim_{k \rightarrow + \infty} \frac{2}{r^k - 1} (1 + q + q^2 + \dots + q^{k - 1}) = 0 .$$
So, given $\epsilon > 0$, we can get $k \in \mathbb{N}^*$ such that
\begin{equation}\label{eq1}
\frac{2}{r^k - 1} (1 + q + q^2 + \dots + q^{k - 1}) < \epsilon .
\end{equation}
Let us fix such a $k$ and put $s := r^k$. According to the theorem 1
(or 2.1) of \cite{far}, the set
$$E_s := \left\{\lambda > 0 ~|~ \norm{\lambda s^n} \leq \frac{1}{s - 1} ~~ (\forall n \in \mathbb{N})\right\}$$
is uncountable. It thus suffices to show that for all $\lambda \in
E_s$, the sequence ${\{<\lambda r^n>\}}_{n \in \mathbb{N}}$ lies in
a finite union of intervals whose sum of lengths $< \epsilon$. Let
$\lambda \in E_s$, arbitrary. We have for all $\ell \in \mathbb{N}$:
$$<\lambda r^{k \ell}> = <\lambda s^{\ell}> \in [0 , \frac{1}{s - 1}] \cup [1 - \frac{1}{s - 1} , 1] .$$
This shows that the sequence ${\{<\lambda r^{k \ell}>\}}_{\ell \in
\mathbb{N}}$ lies in a union of two intervals whose sum of lengths
equal to $\frac{2}{s - 1}$. Now, using the elementary fact stating
that
\begin{quote}
``When the fractional part of a real $x$ lies in a finite union of
intervals whose sum of lengths $\leq \alpha$ $(\alpha
> 0)$ then, given $a , b \in \mathbb{N}^*$, the fractional part of the real $\frac{a}{b} x$ lies in a finite union of intervals
whose sum of lengths $\leq a \alpha$'',
\end{quote}
we deduce that for all $u \in \{1 , 2 , \dots , k - 1\}$, the
sequence ${\{<\lambda r^{k \ell - u}>\}}_{\ell \geq 1}$ lies in a
finite union of intervals whose sum of lengths $\leq \frac{2 q^u}{s
- 1}$. Since each $n \in \mathbb{N}$ can be written in one of the
forms: $k \ell$ $(\ell \in \mathbb{N})$, $k \ell - 1 , k \ell - 2 ,
\dots , k \ell - (k - 1)$ $(\ell \in \mathbb{N}^*)$, it follows that
the sequence ${\{<\lambda r^n>\}}_{n \in \mathbb{N}}$ lies in a
finite union of intervals whose sum of lengths
$$\leq \frac{2}{s - 1} (1 + q + q^2 + \dots + q^{k - 1}) < \epsilon \text{\hspace{1cm}(from (\ref{eq1})),}$$
as required. The proof is complete.
\penalty-20\null\hfill$\blacksquare$\par\medbreak \noindent
 {\bf Remark 1:} If we replace in Theorem \ref{t1} the expression
``uncountable'' by ``infinite'' and we do not require with
$A(\epsilon)$ to have a finite border, the statement which we
obtain is obvious. Indeed, it suffices to pick a
\underline{countable} infinite set $\Gamma \subset ]0 , + \infty[$
and to take:
 $$A(\epsilon) := \{<\lambda r^n> ~|~ \lambda \in \Gamma , n \in \mathbb{N}\} .$$
 Then, we would have $\# A(\epsilon) \leq \# (\Gamma \times
 \mathbb{N})$, which implies that $A(\epsilon)$ is either finite or
 uncountable, consequently it has Lebesgue
 measure zero.~\vspace{1mm}\\
 {\bf Remark 2:} Contrary to the set $\mathbb{Q}_+^*$ of rational positive numbers, the set
 $$\mathbb{Q}_1 := \{p / q ~|~ p , q \in \mathbb{N}^* , p > q^2\}$$
 is representable by an increasing sequence. This follows from the fact that for all
$n \in \mathbb{N}^*$, the set $\mathbb{Q}_1 \cap [1 , n]$ is finite.
In order to prove this for a given $n$, let $p / q \in \mathbb{Q}_1
\cap [1 , n]$, where $p , q \in \mathbb{N}^*$. Then we have $q <
\frac{p}{q} \leq n$ and then $p = \frac{p}{q} q \leq
 n^2$. The positive integers $p$ and $q$ are thus bounded, which implies that the set
 $\mathbb{Q}_1 \cap [1 , n]$ is effectively finite and its cardinal is $\leq n^3$. So, we can ask the question about the value (if it exists) of
 the limit of $\frac{\#(\mathbb{Q}_1 \cap [1 ,
 n])}{n^3}$ $(\leq 1)$ when $n$ tends to infinity.\\
 It is an easy exercise to show that we have for any $n \in
 \mathbb{N}^*$:
 $$\#(\mathbb{Q}_1 \cap [n , n + 1[) = \varphi(1) + \varphi(2) + \dots + \varphi(n)$$
 where $\varphi$ being the Euler indicator. Since we know that $\varphi(1) + \varphi(2) + \dots + \varphi(n)
 = \frac{3 n^2}{\pi^2} + O_n(n \log{n})$ (see e.g. §6 of \cite{hua1}), we deduce that
 $$\#(\mathbb{Q}_1 \cap [1 , n]) \sim_n \frac{1}{\pi^2} n^3 .$$
 The required limit thus exists and equal to $\frac{1}{\pi^2}$. In other words, if $r_n$ denotes the $n$\textsuperscript{th} element of
$\mathbb{Q}_1$ for the increasing order, then we have $r_n \sim
\pi^{2 / 3} \sqrt[3]{n}$. To conclude, this is the first elements
of $\mathbb{Q}_1$:
 $$2 , \frac{5}{2} , 3 , \frac{10}{3} , \frac{7}{2} , \frac{11}{3} , 4 , \frac{17}{4}
 , \frac{13}{3} , \frac{9}{2} , \frac{14}{3} , \frac{19}{4} , 5 , \dots \text{etc} .$$
 {\bf Open problem:} Theorem \ref{t1} is it true for any rational number $p / q > 1$? The case $p / q = 3 / 2$ is particulary
 interesting.

\end{document}